\documentclass{amsart}
\newtheorem{Theorem}{Theorem}[section]
\newtheorem{Definition}[Theorem]{Definition}
\newtheorem{Proposition}[Theorem]{Proposition}

\newtheorem{Lemma}[Theorem]{Lemma}
\newtheorem{Corollary}[Theorem]{Corollary}
\theoremstyle{remark}

\newtheorem{Example}[Theorem]{Example}

\def\il{\int\limits_}
\def\hvi{\varphi}
\def\eps{\varepsilon}

\def\ovr{\overline}

\def\al{\alpha}
\def\gm{\gamma}
\def\th{\theta}

\def\Dl{\Delta}
\def\dl{\delta}

\def\bd{\partial}
\def\lm{\lambda}

\def\sbs{\subset}

\def\dist{\operatorname{dist}}

\def\re{{\mathbf {Re\,}}}

\def\be{\begin{enumerate}}
\def\ee{\end{enumerate}}
\def\bT{\begin{Theorem}}
\def\eT{\end{Theorem}}
\def\bP{\begin{Proposition}}
\def\eP{\end{Proposition}}
\def\bD{\begin{Definition}}
\def\eD{\end{Definition}}
\def\bE{\begin{Example}}
\def\eE{\end{Example}}
\def\bL{\begin{Lemma}}
\def\eL{\end{Lemma}}
\def\bC{\begin{Corollary}}
\def\eC{\end{Corollary}}

\def\oD{\ovr{\mathbb D}}
\def\aD{\mathbb D}
\def\aT{\mathbb T}
\def\aC{\mathbb C}

\def\E{{\mathcal E}}

\begin{document}

\title{On weighted Hardy spaces on the unit disk}
\author{Evgeny A. Poletsky}
\address{Department of Mathematics,  Syracuse University, \newline
215 Carnegie Hall, Syracuse, NY 13244, US\\ Email: eapolets@syr.edu}
\author{Khim R. Shrestha}
\address{Department of Mathematics,  Syracuse University, \newline
215 Carnegie Hall, Syracuse, NY 13244, US, \\ Email: krshrest@syr.edu}

\maketitle
\begin{abstract}
\par In this paper we completely characterize those weighted Hardy spaces that are Poletsky--Stessin Hardy spaces $H^p_u$. We also provide a reduction of $H^\infty$ problems to $H^p_u$ problems and demonstrate how such a reduction can be used to make shortcuts in the proofs of the interpolation theorem and corona problem.
\end{abstract}
\section{Introduction}
\par Let $\lm$ be the normalized Lebesgue measure on the unit circle $\aT$. Among many different definitions of weighted Hardy spaces the closest to our purpose is the definition in \cite{M} and \cite{BPST}. Let $\al\in L^1(\aT)$ be a non-negative function such that $\log\al\in L^1(\aT)$. Then $L^p_\al(\aT)$ is the space of all functions with the finite norm \[\|\phi\|_{\al,p}=\left(\int_0^{2\pi}|\phi(e^{i\th})|^p\al(e^{i\th})\,d\lm\right)^{1/p}\] for $0<p<\infty$ and
$H^p_\al=N^+\cap L^p_\al(\aT)$, where $N^+$ is the Smirnov class. If $\al\equiv1$ then we will use notations $H^p$ and $\|\cdot\|_p$.
\par Later for a plurisubharmonic exhaustion function $u$ on a hyperconvex domain $D\sbs\aC^n$ Poletsky and Stessin introduced in \cite{PS} weighted Hardy spaces $H^p_u(D)$ to generalize the notion of the classical Hardy spaces  and then studied the composition operators generated by the holomorphic mappings between such domains.
\par Recently, M. Alan and N. Gogus in \cite{AG}, S. Sahin in \cite{Sa} and K. R. Shrestha in \cite{Sh1,Sh2} obtained the description of these spaces on the unit disk and called them Poletsky--Stessin Hardy spaces. In particular, they showed that the new spaces form a subclass of weighted Hardy spaces introduced at the beginning of this section.
\par In Section \ref{PSH} we finish this description proving that $H^p_\al$ is a Poletsky--Stessin space if and only if the weight $\al$ is lower semicontinuous and greater than some $c>0$ on $\aT$.
\par Although weighted Hardy spaces can be studied per se there is also an expectation that they can be useful for the classical theory. In Section \ref{S:hp} we prove the main result of this paper: if a closed convex set $A$ intersects unit balls in all $H^p_u(\aD)$ for some $p>1$ then it intersects the unit ball in $H^\infty$. Thus to find bounded solutions to a linear problem it suffices to show that they exist at all $H^p_u(\aD)$ and their norms are uniformly bounded.
\par In the last two section we use this fact to demonstrate shortcuts in the proofs of the interpolation theorem and corona problem.
\section{Duality}\label{S:WHS}
\par
\par Let $a(z)$ be a holomorphic function such that $|a(e^{i\th})|=\al(e^{i\th})$ on $[0,2\pi]$ a.e. and $a$ never takes the zero value. Such function does exist and belongs to $H^1$ because $\log\al$ is integrable on $\aT$ so we can take a harmonic function
\[h(z)=\frac1{2\pi}\il0^{2\pi}\log\al(e^{i\th})P(z,e^{i\th})\,d\th,\]
add a conjugate function $g$ and write $a(z)=e^{h(z)+ig(z)}$.
\par In \cite{BPST} for $f\in L^p_\al(\aT)$ the operator $A_pf=a^{1/p}f$ was introduced. Then
\[\|A_pf\|_p^p=\il0^{2\pi}|f(e^{i\th})|^p\al(e^{i\th})\,d\lm=
\|f\|^p_{\al,p}.\]
Thus $A_p$ is an isometrical imbedding of $L^p_\al(\aT))$ into $L^p(\aT)$.
\par We add another requirements on the weight $\al$ asking that $\al\ge c>0$ on $\aT$.
Clearly, $A_p^{-1}f=\al^{-1/p}f$ is also an isometry and the inverse of $A_p$. Hence $A_p$ is an isometric isomorphism of $L^p_\al(\aT))$ onto $L^p(\aT)$.  Moreover, $A_p$ maps $H^p_\al$ isometrically onto $H^p$.
\par If $\phi\in L^p_\al(\aT)$ then $\dist(\phi,H^p_\al)=\dist(A_p\phi,H^p)$. By the classical result (see \cite{K}) for $p>1$
\begin{equation}\label{e:ds}
\dist(A_p\phi,H^p)=\left|\sup_{g\in H^q,\|g\|_{H^q}=1}
\il0^{2\pi}a^{1/p}(e^{i\th})\phi(e^{i\th})e^{i\th}g(e^{i\th})\,d\lm\right|.
\end{equation}
\par Since the $H^q_\al$-norm of $\al^{-1/q}g(z)$ coincides with the $H^q$-norm of $g$ we can get the following duality result:
\bT\label{T:dual} If $\phi\in L^p_\al(\aT)$ then
\begin{align*}
\dist(\phi,H^p_\al)=\left|\sup_{g\in H^q_\al,\|g\|_{H^q_\al}=1}\il{\aT}\phi(e^{i\th})a(e^{i\th})e^{i\th}g(e^{i\th})\,d\lm\right|.
\end{align*}
\eT
\section{Poletsky--Stessin Hardy spaces} \label{PSH}

\par Let $\aD$ be the unit disc $\{|z|< 1\}$ in $\aC$.  A continuous subharmonic function $u:\mathbb{D}\to [-\infty, 0)$ such that $u(z)\to 0$ as $|z|\to 1$ is called an exhaustion function. The class of such functions will be denoted by $\E$. Following \cite{D} for $r <0$ we set
\[B_{u,r}  =\{z\in \mathbb{D}:u(z) < r\}\text{ and }
S_{u,r} = \{z\in \mathbb{D}: u(z) = r\}.\]
As in \cite{D} we let $u_r = \max\{u,r\}$ and define the  measure
\[\mu_{u,r} = \Dl u_r -\chi_{\aD \setminus B_r}\Dl u, \] where $\Dl$ is the Laplace operator.
Clearly $\mu_{u,r}\ge 0$ and is  supported by $S_{u,r}$.
\par Let us denote by $\E_1$ the set of all continuous negative subharmonic exhaustion functions $u$ on $\aD$ such that
\[\il \aD\Dl u = 1.\] In the same paper Demailly (see Theorems 1.7 and 3.1 there) proved the following result which we adapt to the case of $\aD$.
\bT [Lelong--Jensen formula]\label{jlf}
Let $\phi$ be a subharmonic function on $\mathbb{D}$. Then $\phi$ is $\mu_{u,r}$-integrable for every $r < 0$ and \[\int_{S_{u,r}}\phi\,d\mu_{u,r}= \int_{B_{u,r}} \phi\, \Dl u + \int_{B_{u,r}}(r-u)\,\Dl\phi.\]
Moreover, if $u\in\E_1$ then the measures $\mu_{u,r}$ converge weak-$*$ in $C^*(\oD)$  to a measure $\mu_u\ge 0$ supported by $\aT$ as $r\to 0^-$.
\eT
As a consequence of this theorem we have the following corollary.
\bC\label{C:imu} If $\phi$ is a non-negative subharmonic function, then the function
\[r\mapsto\int_{S_{u,r}}\phi\,d\mu_{u,r}\] is increasing on $(-\infty, 0)$.\eC
\par Using the measures $\mu_{u,r}$ Poletsky and Stessin introduced in \cite{PS} the weighted Hardy spaces associated with an exhaustion $u\in\E$. Following \cite{PS} we define the space $H^p_u, 0 < p < \infty,$ consisting of all holomorphic functions $f(z)$ in $\mathbb{D}$ that satisfy \[\|f\|_{u,p}^p = \varlimsup_{r\to 0^-}\il {S_{u,r}} |f|^p\, d\mu_{u,r} < \infty.\] By Corollary \ref{C:imu} we can replace the $\varlimsup$ in the above definition with $\lim$.  By Theorem \ref{jlf} and the monotone convergence theorem it follows that,
\begin{equation}\label{e:nf} \|f\|^p_{u,p} = \il {\aD} |f|^p\,\Dl u - \il {\aD} u\,\Dl |f|^p.
\end{equation}
The classical Hardy spaces $H^p$ correspond to the exhaustion function $u(z)=\log|z|$ (\cite[Section 4]{PS}). Hence the classical definition of the Hardy spaces is subsumed in this new definition.
\par It is proved in \cite{PS} that the spaces $H^p_u$ are Banach when $p\ge 1$ and
if $v,u\in\E$ and $v\le cu$ in a neighborhood of $\aT$ for some $c>0$, then $H^p_v\sbs H^p_u$ and if $f\in H^p_v$ then $\|f\|_{u,p}^p\le c\|f\|_{v,p}^p$. Thus by Hopf's lemma the space $H^p_u(\mathbb{D})$ is contained in the classical Hardy space $H^p(\mathbb{D})$.
\par  It has been established (see \cite{AG}, \cite{Sa} and \cite{Sh1}) that the boundary measure $\mu_u=\al_u d\lm$ for some $\al_u \in L^1(\aT)$ and $f\in H^p_u$ if and only if $f\in H^p$ and $\|f\|_{p,\al_u}<\infty$. Moreover, $\|f\|_{u,p}=\|f\|_{\al_u,p}$. Hence
\begin{equation}\label{e:nf1} \il{\aT}|f|^p\,d\mu_u=\il {\aD} |f|^p\,\Dl u - \il {\aD} u\,\Dl |f|^p.
\end{equation}
\par The weight $\al_u$ has the following properties:
\begin{enumerate}
\item \[\al_u(e^{i\th})=\il {\aD} P(z, e^{i\th})\,\Dl u(z),\] where $P(z, e^{i\th})$ is the Poisson kernel;
\item $\|\al_u\|_{L^1} = 1$ if and only if $u\in\E_1$;
\item $\al_u(e^{i\th})$ is lower semicontinuous and $\al_u(e^{i\th})\ge c $ on $\aT$ for some $c >0$.
\end{enumerate}
\par The class of Poletsky--Stessin Hardy spaces is more narrow than weighted spaces discussed in Section \ref{S:WHS} just because the weight function $\al_u$ must be lower semicontinuous and be greater than some $c>0$ on $\aT$. As the following result shows these are the only restrictions on weights.
\bT\label{T:wc} Let $\al$ be a measurable function on $\aT$. Then $\al\,d\lm=\mu_u$ for some $u\in\E_1$ if and only if $\al$ is lower semicontinuous, $\al(e^{i\th})\ge c>0$ on $\aT$ and
\begin{equation}\label{e:tww}
\il\aT\al\,d\lm=1.
\end{equation}
\eT
\begin{proof}
\par Let $\al\in C(\aT)$ be a function such that $\al\ge c>0$ on $\aT$. For $0 < r < 1$ define
\[\al_r(e^{i\th}) = \int_\aT P(re^{i\th}, e^{i\hvi})\al(e^{i\hvi})\,d\lm(\hvi).\]
Then $\al_r \to \al$ uniformly on $\aT$ as $r\to 1$. Clearly $\al_r \in C^\infty(\aT)$.
\par Define \[ u_r(z) = \int \log\left|\frac{z - re^{i\hvi}}{1 -re^{-i\hvi}z}\right|\al(e^{i\hvi})\,d\lm(\hvi).\]
Then $u_r$ is a subharmonic exhaustion function on $\aD$ and by the Riesz Decomposition Theorem its Laplacian $\Dl u_r$ is supported by $\aT(r)=\{z=re^{i\phi}\}$ and is equal to
$\al(e^{i\hvi})\,d\lm(\hvi)$.
Hence
\[\int_\aD \Dl u_r(z) = \int_\aT\al(e^{i\hvi})\,d\lm(\hvi).\]
\par The weight $\al_r(e^{i\th})$ of $u_r$ is equal to
\[\il {\aT}P(re^{i\hvi}, e^{i\th})\al(e^{i\hvi})\,d\lm(\hvi)=\al_r(e^{i\th}).\]
Hence any $\al\in C(\aT)$ can be uniformly approximated by a function $\beta_u$ such that $\beta_u\,d\lm=\mu_u$ and $u\in\E$.
\par If $\al$ is any lower semicontinuous function satisfying (\ref{e:tww})  and such that $\al\ge c>0$ on $\aT$, then $\al$ is the pointwise limit of an increasing sequence of continuous functions $\al_j$ such $\al_j\ge c/2>0$ on $\aT$. Replacing $\al_j$ with the functions $\al_j-2^{-j}$ we may assume that the function $\beta_j=\al_{j}-\al_{j-1}\ge 2^{-j-1}$ on $\aT$. (Here we set $\al_0=0$.) By the argument above we can approximate the functions $\beta_j$ by continuous functions $\gm_j$ such that $\gm_j\ge2^{-j-2}$ on $\aT$, $\gm_j\,d\lm=\mu_{u_j}$ for some $u_j\in\E$ and
\[\sum_{j=1}^\infty\gm_j=\al.\]
\par Let $v_j=\max\{u_j,-2^{-j}\}$. Since for a fixed $j$ the weak-$*$ limits of $\mu_{u_j,r}$ and $\mu_{v_j,r}$ as $r\to0^-$ coincide we see that $\al_{v_j}=\al_{u_j}=\gm_j$. If $v=\sum v_j$ then $v$ is a continuous exhaustion of $\aD$ such that $\lim_{|z|\to1}v(z)=0$. Moreover,
\[\il{\aD}\Dl v=\sum_{j=1}^\infty\il{\aD}\Dl v_j=\sum_{j=1}^\infty\il{\aT}\gm_j=\il{\aT}\al=1.\]
Hence $v\in\E_1$.
\par Now
\[\il{\aD}P(z,e^{i\th})\Dl v(z)=\sum_{j=1}^\infty\il{\aD}P(z,e^{i\th})\Dl v_j(z)=\sum_{j=1}^\infty\gm_j(e^{i\th})=\al(e^{i\th}).\]
Thus $\mu_v=\al$.
\par The converse statements follows from results \cite{AG}, \cite{Sa} and \cite{Sh1} mentioned above.
\end{proof}
\par Among the advantages of these spaces comparatively to spaces studied in \cite{BPST} we can list the following. First of all, one does not need the existence of boundary values or the notion of Smirnov class to introduce these spaces. This is especially attractive for the theory of functions in several variables on non-smooth domains.
\par Another advantage is the existence of Carleson measures.  Given a weight $\al$ a measure $\nu$ on the unit disk $\aD$ is called $\al$-Carleson with the constant $C(\al)$ if \[\il{\aD}|f|^p\,d\nu\le C(\al)\il{\aT}|f|^p\al\,d\lm\] for all $p>1$ and all $f\in H^p_\al$. If $\al\equiv1$ then such measure are called Carleson measures. In \cite{M} one can find the characterisation of $\al$-Carleson measures for $\al$ satisfying Muckenhoupt's conditions similar to the classical characterisation of Carleson measures by L. Carleson. In the case of Poletsky--Stessin Hardy spaces it follows immediately from (\ref{e:nf1}) that the measure $\Dl u$ is $\al_u$-Carleson with the constant $C(\al_u)=1$. By Theorem \ref{T:wc} we see that $\al$-Carleson measures with constant 1 exist for all lower semicontinuous weights.
\par Thirdly, the formula (\ref{e:nf1}) helps to obtain additional information. For example, one can get integrability of derivative. Since $\Dl|f|^p=\frac{p^2}4|f|^{p-2}|f'|^2$ for all $f\in H^p_u$, $p\ge 1$, we have the inequality
\[\il{\aT}|f|^p\,d\mu_u\ge\frac{p^2}4\il {\aD}|u||f|^{p-2}|f'|^2\,dx\,dy.\]
\section{From $H^p_u$ to $H^\infty$}\label{S:hp}
\par Let $u_1,\dots,u_k$ be exhaustion functions from $\E_1$ and let $u=(u_1,\dots,u_k)$. We say that $u\in\E^k_1$. Let $H^p_u$ to be the direct product $H^p_{u_1}\times\cdots\times H^p_{u_k}$ with the norm
\[\|(f_1,\dots,f_k)\|_{u,p}=\sum_{j=1}^k\|f_j\|_{u_j,p}.\] We will use the notation $(H^p)^k$ and $\|f\|_p$ when $\al_{u_1}=\cdots\al_{u_k}=1$. We denote by $B_{u,p}(r)$ the closed ball of radius $r$ centered at the origin of $H^p_u$
\par  The norm on $(H^\infty)^k$ will be defined as
\[\|f\|_{\infty}=\sum_{j=1}^k\|f_j\|_\infty\] and $B_\infty(r)$ is the closed ball of radius $r$ centered at the origin of $(H^\infty)^k$. Then $B_\infty(r)\sbs B_{u,p}(r)$.
\bT \label{T:int1}Let $A \sbs (H^p)^k$, $p>1$, be a closed convex set. Then $A\cap B_\infty(1) \neq \emptyset$ if and only if $A\cap B_{u,p}(1) \neq \emptyset$ for all exhaustion vector-functions $u=(u_1,\dots,u_k)\in\E^k_1$. \eT
\begin{proof} The ``only if" part of the theorem is obvious. The ``if'' part will be proved by contradiction. Let us take $0<\eps<1$ and suppose that $A\cap B_\infty (r_0)=\emptyset$ for $r_0=(1-\eps)^{-1}$. By the Hahn--Banach theorem there exists $g=(g_1,\dots,g_k)\in (L^q(\aT))^k$ such that
\[\sum_{j=1}^k\re\il{\aT} f_jg_j\,d\lm \ge 1\] for all $f\in A$ and
\[ \sum_{j=1}^k\re\il{\aT} f_jg_j\,d\lm \le 1\] for all $f\in B_\infty(r_0)$. Multiplying $f_j$ by appropriate constants $a_j$ with $|a_j|=1$ we see that
\[\sum_{j=1}^k\left|\il{\aT} f_jg_j\,d\lm\right| \le r^{-1}_0=1-\eps\] for all $f\in B_\infty(1)$.
\par Let $\tilde{g}_j(z) = g_j(z)/z$.  Then $\tilde{g}_j\in L^q(\aT)\sbs L^1(\aT)$ for all $j$.  By a duality result (see \cite[VII.2]{K}) there exist $h_j\in H^1$ and $p_j\in H^\infty$ such that $\|p_j\|_\infty = 1$, $p_j(0) = 0$ and \[(\tilde{g}_j - h_j)p_j = |\tilde{g}_j-h_j|\]
almost everywhere.
\par We take $f=(f_1,\dots,f_k)\in (H^\infty)^k$ such that $f_i\equiv0$ when $i\ne j$ and $f_j(z)=p_j(z)/z$. Clearly, $f\in B_\infty(1)$. Therefore,
\[1-\eps\ge\left|\il{\aT} f_jg_j\,d\lm\right|=\left|\il{\aT} (\tilde{g}_j-h_j)p_j\,d\lm\right|=\il{\aT}|\tilde{g}_j-h_j|\,d\lm.\]
\par There is $\tilde h_j\in H^q$  so that $\|h_j-\tilde h_j\|_1 \le\eps/2$.  Let $\phi_j=|\tilde{g}_j-\tilde h_j|$. Then
\[\il{\aT}\phi_j\,d\lm\le\il{\aT}\left(|\tilde{g}_j-h_j|+|h_j-\tilde h_j|\right)\,d\lm\le 1-\eps/2.\]
\par And for $f\in A$,
\[\sum_{j=1}^k\il{\aT}\phi_j|f_j|\,d\lm=\sum_{j=1}^k\il{\aT}|(g_j-z\tilde h_j)f_j|\,d\lm \ge \sum_{j=1}^k\left|\il{\aT}(g_j-z\tilde h_j)f_j\,d\lm\right| \ge 1.\]
\par Let $\tilde\phi_j=\max\{\phi_j,\eps/4\}$.
Then $\|\tilde\phi_j\|_1\le\|\phi_j+\eps/4\|_1\le 1-\eps/4$. Now for $f\in A$,
\[\sum_{j=1}^k\il{\aT}|f_j|\tilde\phi_j\,d\lm\ge\sum_{j=1}^k\il{\aT}|f_j|\phi_j\,d\lm\ge1.\]
\par For any $\dl>0$ and $1\le j\le k$ there exists $\psi_j\in C(\aT)$ such that  $\psi_j\ge\eps/8$, $\|\psi_j\|_1 =\|\tilde\phi_j\|_1$ and $\|\psi_j-\tilde\phi_j\|_q<\dl$. \par For $f\in A$,
\[\sum_{j=1}^k\il{\aT}|f_j|\psi_j\,d\lm\ge\sum_{j=1}^k\il{\aT}|f_j|\tilde\phi_j\,d\lm-
\sum_{j=1}^k\il{\aT}|f_j||\psi_j- \tilde\phi_j|\,d\lm\ge 1-\dl\|f\|_p.\]
\par By Theorem \ref{T:wc} there are exhaustion functions $u_j$, $1\le j\le k$, such that $\mu_{u_j}=a_j\psi_j$, where $a_j$ is chosen so that $\|a_j\psi_j\|_1=1$. Let $u=(u_1,\dots,u_k)$. Note that $a_j\ge(1-\eps/4)^{-1}$.
\par If $f\in B_{u, p}(1)$ then
\[\sum_{j=1}^k\il{\aT}|f_j|a_j\psi_j\,d\lm\le
\sum_{j=1}^k\|f_j\|_{u,p}\|a_j\psi_j\|_1^{1/q}\le1\]
and
\[\|f\|_p=\sum_{j=1}^k\left(\il{\aT}|f_j|^p\,d\lm\right)^{1/p} \le \left(\frac{\eps}{8}\right)^{-1/p}\sum_{j=1}^k\|f_j\|_{u_j,p}\le
\left(\frac{\eps}{8}\right)^{-1/p}=c.\]
Thus if $f\in A$ and $\|f\|_p>c$ then $f\not\in B_{u, p}(1)$. On the other hand if $f\in A$ and $\|f\|_p \le c$, then
\[\sum_{j=1}^k\il{\aT}|f_j|a_j\psi_j\,d\lm\ge (1-\eps/4)^{-1}(1-c\dl).\]
\par Taking $\dl>0$ so small that $(1-\eps/4)^{-1}(1-c\dl)>1$ we see that $A\cap B_{u,p}(1)=\emptyset$. Hence $A\cap B_\infty(r_0)\ne\emptyset$ for all $r_0>1$.
\par Let $\{f_n\}$ be a sequence of functions such that $f_n\in A\cap B_\infty(1+1/n)$. We may assume that $\{f_n\}$ converges uniformly on compacta to a function $f\in B_\infty(1)$. This implies that $\{f_n\}$ converges to $f$ weakly. Since any convex closed set is weakly closed we see that $f\in A$.
\end{proof}
\par As the following corollary shows it is possible to use the theorem above when all functions $u_j$ are equal although constants will change.
\bC \label{C:int1}Let $A \sbs (H^p)^k$, $p>1$, be a closed convex set. Suppose $A\cap B_{{\bf u},p}(1) \neq \emptyset$ for all exhaustion vector-functions ${\bf u}=(u,\dots,u)\in\E^k_1$. Then $A\cap B_\infty(k)\neq\emptyset$. Conversely, if $A\cap B_\infty(1) \neq \emptyset$ then $A\cap B_{{\bf u},p}(1)\neq\emptyset$.\eC
\begin{proof} Let $v=(v_1,\dots,v_k)\in\E_1^k$. Let
\[u=\frac1k\sum_{j=1}^kv_j.\]
Then $u\in\E_1$ and by the assumption of the corollary there is $f=(f_1,\dots,f_k)\in A\cap B_{{\bf u},p}(1)$, where ${\bf u}=(u,\dots,u)$. Note that $v_j\ge ku$. By Corollary 3.2 in \cite{PS} $\|f_j\|_{v_j,p}\le k\|f_j\|_{u,p}$, $1\le j\le k$. Hence $f\in B_{v,p}(k)$ and $A\cap B_{v,p}(k) \neq \emptyset$. By Theorem \ref{T:int1} $A\cap B_\infty(k)\neq\emptyset$.
\end{proof}
\section{Interpolation Theorem}\label{S:it}
\par A sequence $\{z_j\}_1^\infty\sbs\aD$ is $\dl$-sparse for $\dl>0$ if
\[\inf_{k}\prod_{j\ne k}|\frac{z_j-z_k}{1-\bar z_kz_j}|\ge\dl\] for all $k$.
\par A sequence $\{z_j\}\sbs\aD$ is called interpolating if for any sequence $s=\{s_j\}\in l^\infty$ there is a function $f\in H^\infty$ such that $f(z_j)=s_j$ for all $j$ and $\|f\|_{H^\infty}\le C\|s\|_\infty$ and the constant $C$ does not depend on $\|s\|_\infty$.
\par The famous theorem of Carleson states
\bT\label{T:it} A sequence $\{z_j\}\sbs\aD$ is interpolating if and only if it is $\dl$-sparse for some $\dl>0$.
\eT
\par Now we can present a shorter proof of the sufficiency part of Theorem \ref{T:it}. The proof of necessity is quite elementary and can be found in \cite{G}.  Theorem 3.2 in \cite{H} that is a quick consequence of the general characterization of Carleson measures, states that if a sequence $\{z_j\}\sbs\aD$ is $\dl$-sparse then  the measure \[\nu=\sum_{j=1}^\infty(1-|z_j|^2)\dl_{z_j}\] is Carleson with a constant $C$ depending only on $\dl$.
\par We take an integer $N>1$ and denote by $X_N$ the set of all functions  $f\in H^2$ such that $f(z_j)=s_j$, $1\le j\le N$. Clearly $X_N$ is closed and convex.
\par Let
\[B(z)=\prod_{j=1}^N\frac{z-z_j}{1-\bar z_jz}\text{ and } B_k(z)=\prod_{j=1,j\ne k}^N\frac{z-z_j}{1-\bar z_jz},\quad k=1,\dots,N.\]
Then any function $f$ in $X_N$ has the form
\[\sum_{j=1}^N\frac{s_j}{B_j(z_j)}B_j(z)+B(z)h(z)=
\left(\sum_{j=1}^N\frac{s_j}{B_j(z_j)}\frac{1-\bar z_jz}{z-z_j}+h(z)\right)B(z),\]
where $h\in H^2$.
\par We set $C_j=s_jB^{-1}_j(z_j)$ and let
\[\phi(z)=\sum_{j=1}^NC_j\frac{1-\bar z_jz}{z-z_j}.\]
Let $u\in\E_1$ and let $a=a_u$ be the function introduced in Section \ref{S:WHS}. Then for $g\in H^2$ with $\|g\|_{H^2}=1$ we have
\begin{align*}
&\frac1{2\pi}\il{\aT}\phi(z)a^{1/2}(z)g(z))\,dz=
\sum_{j=1}^NC_j(1-|z_j|^2)g(z_j)a^{1/2}(z_j)\\&\le
\frac{\|s\|_\infty}{\dl}\il{\aD}|ga^{1/2}|\,d\nu\le\frac{\|s\|_\infty}{\dl}
\left(\il{\aD}|g|^2\,d\nu\right)^{1/2}\left(\il{\aD}|a|\,d\nu\right)^{1/2}
\\&\le\frac{C^2\|s\|_\infty}{\dl}\|g\|_{H^2}\|a^{1/2}\|_{H^2}=
\frac{C^2\|s\|_\infty}{\dl}=C'\|s\|_\infty.
\end{align*}
Hence by (\ref{e:ds}) $\dist(\phi,H^2_u)\le C'\|s\|_\infty$ and this means that $X_N\cap B_{u,2}(C'\|s\|_\infty)\ne\emptyset$. Thus by Theorem \ref{T:int1} there is $f_N\in X_N\cap B_\infty(C'\|s\|_\infty)$. Since $C'$ does not depend on $N$ there is $f\in B_\infty(C'\|s\|_\infty)$ interpolating $s$.
.
\section{Corona Theorem}\label{S:ct}
\par The same method can be applied to the corona theorem.
\bT If the functions $f_1,\dots,f_n$ are in the unit ball of $H^\infty$ and
\[\sum_{j=1}^n|f_j|^2\ge\dl>0,\]then there are functions $g_1,\dots,g_n$ in $H^\infty$ such that
\begin{equation}\label{e:ce}
\sum_{j=1}^nf_jg_j=1
\end{equation} and $\|g_j\|\le C$, where $C$ depends only on $\dl$.
\eT
\par We will discuss only the case when $n=2$. It suffices to prove this theorem for  functions $f_j$ that can be continuously extended to $\oD$ and have finitely many zeros in $\oD$. In this case one can easily find functions $\phi_1$ and $\phi_2$ smooth up to the boundary such that
\[f_1\phi_1+f_2\phi_2=1.\]
To make them holomorphic we look for a function $v$ such that
\[\bar\bd(\phi_1+f_2v)=\bar\bd\phi+f_2\bar\bd v=0\]
and
\[\bar\bd(\phi_2-f_1v)=\bar\bd\phi-f_1\bar\bd v=0.\]
Since $f_1\bar\bd\phi_1+f_2\bar\bd\phi_2=0$ we see that
\[\bar\bd v=f_1^{-1}\bar\bd\phi_2=-f_2^{-1}\bar\bd\phi_1=:\psi.\]
\par The following lemma can be found in \cite{G}.
\bL There are solutions $\phi_1$ and $\phi_2$ to (\ref{e:ce}) continuous up to the boundary such that the measure $\nu=|\psi|\,dzd\bar z$ is Carleson with constant $C$ depending only on $\dl$ and $|\phi_1|+|\phi_2|\le K(\dl)$.
\eL
\par Let
\[\Psi(z)=\il{\aD}\frac{\psi(\zeta)}{\zeta-z}\,d\zeta d\bar\zeta.\]
Then $\bar\bd\Psi=\psi$ and  for any $u\in\E_1$
\begin{align*}&\left|\il{\aT}\Psi(z)a_u(z)g(z)\,dz\right|^2=
\left|\il{\aD}\psi(\zeta)a_u(\zeta)g(\zeta)\,d\zeta d\bar\zeta\right|^2\\&\le
\il{\aD}|a_u(\zeta)|\psi(\zeta)\,d\zeta d\bar\zeta \il{\aD}|\psi(\zeta)||a_u(\zeta)g^2(\zeta)|\,d\zeta d\bar\zeta\le C^2\|g\|^2_{u,2}.
\end{align*}

Thus by Theorem \ref{T:dual} $\dist(\Psi,H^2_u)\le C$.
Hence there is $v=\Psi+h$ such that $h\in H^2_u$ and $\|h\|_{H^2_u}\le C$. Therefore the function $h_1=\phi_1+f_2v$ is holomorphic, lies in $H^2_u$ and $\|h_1\|_{u,2}\le K(\dl)+C=R$. The same estimate holds for the function $h_2=\phi_2-f_1v$.
\par Thus if $A\in (H^2)^2$ is the set of all solutions $(g_1,g_2)$ to (\ref{e:ce}), then $A\cap B_{u,2}(R)\ne\emptyset$ for all pairs $(u,u)$, where $u\in\E_1$. Since the set $A$ is convex and closed, by Corollary \ref{C:int1} $A\cap B_\infty(2R)\ne\emptyset$. This ends the proof.

\end{document}